\newfont{\bcb}{msbm10}
\newfont{\matb}{cmbx10}
\newfont{\got}{eufm10}

\documentclass[12pt]{amsart}
\usepackage{amsmath, amsthm, amscd, amsfonts, amssymb, latexsym, graphicx, color}
\usepackage[bookmarksnumbered, colorlinks, plainpages, hypertex]{hyperref}

\usepackage[cp1250]{inputenc}

\usepackage{amsmath,amsthm}
\usepackage{amssymb,latexsym}
\usepackage{enumerate}

\newtheorem{theorem}{Theorem}[section]

\newtheorem{proposition}[theorem]{Proposition}
\newtheorem{corollary}[theorem]{Corollary}
\theoremstyle{definition}

\theoremstyle{remark}

\numberwithin{equation}{section}

\begin{document}

\title[Piecewise continuity]{Piecewise continuity
       \\ of functions definable over \\ Henselian rank one valued fields}

\author[Krzysztof Jan Nowak]{Krzysztof Jan Nowak}


\subjclass[2000]{12J25, 13F30, 14P10}

\keywords{}



\begin{abstract}
Consider a Henselian rank one valued field $K$ of
equi\-characteristic zero along with the language
$\mathcal{L}^{P}$ of Denef--Pas. Let $f: A \to K$ be an
$\mathcal{L}^{P}$-definable (with parameters) function on a subset
$A$ of $K^{n}$. We prove that $f$ is piecewise continuous; more
precisely, there is a finite partition of $A$ into
$\mathcal{L}^{P}$-definable locally closed subsets
$A_{1},\ldots,A_{s}$ of $K^{n}$ such that the restriction of $f$
to each $A_{i}$ is a continuous function.
\end{abstract}

\maketitle


\section{Introduction}

Consider a Henselian rank one valued field $K$ of
equi\-characteristic zero along with the language
$\mathcal{L}^{P}$ of Denef--Pas, which consists of three sorts:
the valued field $K$-sort, the value group $\Gamma$-sort and the
residue field $\Bbbk$-sort. The only symbols of $\mathcal{L}^{P}$
connecting the sorts are the following two maps from the main
$K$-sort to the auxiliary $\Gamma$-sort and $\Bbbk$-sort: the
valuation map $v$ and an angular component map $\overline{ac}$
which is multiplicative, sends $0$ to $0$ and coincides with the
residue map on units of the valuation ring $R$ of $K$. The
language of the $K$-sort is the language of rings; that of the
$\Gamma$-sort is any augmentation of the language of ordered
abelian groups (with $\infty$); finally, that of the $\Bbbk$-sort
is any augmentation of the language of rings. Throughout the paper
the word "definable" means "definable with parameters". We
consider $K^{n}$ with the product topology, called the
$K$-topology on $K^{n}$. Let $\mathbb{P}^{1}(K)$ stand for the
projective line over $K$.

\vspace{1ex}

The main purpose is to prove the following

\begin{theorem}\label{piece}
Let $A \subset K^{n}$ and $f: A \to \mathbb{P}^{1}(K)$ be an
$\mathcal{L}^{P}$-definable function in the three-sorted language
of Denef--Pas. Then $f$ is piecewise continuous, i.e.\ there is a
finite partition of $A$ into $\mathcal{L}^{P}$-definable locally
closed subsets $A_{1},\ldots,A_{s}$ of $K^{n}$ such that the
restriction of $f$ to each $A_{i}$ is continuous.
\end{theorem}

We immadiately obtain

\begin{corollary}
The conclusion of the above theorem holds for any
$\mathcal{L}^{P}$-definable function $f: A \to K$.
\end{corollary}

The proof of Theorem~\ref{piece} relies on two basic ingredients.
The first one is concerned with a theory of algebraic dimension
and decomposition of definable sets into a finite union of locally
closed definable subsets, recalled in the next section. It was
established by van den Dries~\cite{Dries} for certain expansions
of rings (and Henselian valued fields, in particular) which admit
quantifier elimination and are equipped with a topological system.
The second one is the closedness theorem from our
paper~\cite[Theorem~3.1]{Now}. Let us mention that the latter
paper is devoted to geometry over Henselian rank one valued fields
and i.a.\ to the results achieved in our joint article~\cite{K-N}
about hereditarily rational functions on real and $p$-adic
varieties. Section~3 provides the proof of the main result
(Theorem~\ref{piece}).

\section{Definable sets over Henselian valued fields}

Consider an infinite integral domain $D$ with quotient field $K$.
One of the fundamental concept introduced by van den
Dries~\cite{Dries} is that of a \emph{topological system} on a
given expansion $\mathcal{D}$ of a domain $D$ in a language
$\mathcal{L}$. That concept incorporates both Zariski-type and
definable topologies. We remind the reader that it consists of a
topology $\tau_{n}$ on each set $D^{n}$, $n \in \mathbb{N}$, such
that:

1) For any $n$-ary $\mathcal{L}_{D}$-terms $t_{1},\ldots,t_{s}$,
$n,s \in \mathbb{N}$, the induced map
$$ D^{n} \ni a \longrightarrow (t_{1}(a),\ldots,t_{s}(a)) \in D^{s}
$$
is continuous.

2) Every singleton $\{ a \}$, $a \in D$, is a closed subset of
$D$.

3) For any $n$-ary relation symbol $R$ of the language
$\mathcal{L}$ and any sequence $1 \leq i_{1} < \ldots < i_{k} \leq
n$, $1 \leq k \leq n$, the two sets
$$ \{ (a_{i_{1}},\ldots,a_{i_{k}}) \in D^{k}: \ \mathcal{D} \models
   R((a_{i_{1}},\ldots,a_{i_{k}})^{\&}), \, a_{i_{1}} \neq 0, \ldots,
   a_{i_{k}} \neq 0 \}, $$
$$ \{ (a_{i_{1}},\ldots,a_{i_{k}}) \in D^{k}: \ \mathcal{D} \models
   \neg R((a_{i_{1}},\ldots,a_{i_{k}})^{\&}), \, a_{i_{1}} \neq 0, \ldots,
   a_{i_{k}} \neq 0 \} $$
are open in $D^{k}$; here $(a_{i_{1}},\ldots,a_{i_{k}})^{\&}$
denotes the element of $D^{n}$ whose $i_{j}$-th coordinate is
$a_{i_{j}}$, $j=1,\ldots,k$, and whose remaining coordinates are
zero.

\vspace{1ex}

Finite intersections of closed sets of the form
$$ \{ a \in D^{n}: t(a)=0 \}, $$
where $t$ is an $n$-ary $\mathcal{L}_{D}$-term, will be called
\emph{special closed subsets} of $D^{n}$. Finite intersections of
open sets of the form
$$ \{ a \in D^{n}: t(a) \neq 0 \}, $$
$$ \{ a \in D^{n}: \ \mathcal{D} \models
   R((t_{i_{1}}(a),\ldots,t_{i_{k}}(a))^{\&}), \, t_{i_{1}}(a) \neq 0, \ldots,
   t_{i_{k}}(a) \neq 0 \} $$
or
$$ \{ a \in D^{n}: \ \mathcal{D} \models
   \neg R((t_{i_{1}}(a),\ldots,t_{i_{k}}(a))^{\&}), \, t_{i_{1}}(a) \neq 0, \ldots,
   t_{i_{k}}(a) \neq 0 \}, $$
where $t,t_{i_{1}},t_{i_{k}}$ are $\mathcal{L}_{D}$-terms, will be
called \emph{special open subsets} of $D^{n}$. Finally, an
intersection of a special open and a special closed subsets of
$D^{n}$ will be called a \emph{special locally closed} subset of
$D^{n}$. Every quantifier-free $\mathcal{L}$-definable set is a
finite union of special locally closed sets.

\vspace{1ex}

Suppose now that the language $\mathcal{L}$ extends the language
of rings and has no extra function symbols of arity $>0$ and that
an $\mathcal{L}$-expansion $\mathcal{D}$ of the domain $D$ under
study admits quantifier elimination and is equipped with a
topological system such that every non-empty special open subset
of $D$ is infinite. These conditions ensure that $\mathcal{D}$ is
algebraically bounded and algebraic dimension defines a dimension
function on $\mathcal{D}$ (\cite[Proposition~2.15
and~2.7]{Dries}). Algebraic dimension is the only dimension
function on $\mathcal{D}$ whenever, in addition, $D$ is a
non-trivially valued field and the topology $\tau_{1}$ is induced
by its valuation. Then, for simplicity, the algebraic dimension of
an $\mathcal{L}$-definable set $E$ will be denoted by $\dim E$.

\vspace{1ex}

Now we recall the following two basic results
from~\cite[Proposition~2.17 and~2.23]{Dries}:

\begin{proposition}\label{partition}
Every $\mathcal{L}$-definable set is a finite union of
intersections of Zariski closed with special open sets and, a
fortiori, a finite union of locally closed sets.
\end{proposition}

\begin{proposition}\label{front-eq}
Let $E$ be an $\mathcal{L}$-definable subset of $D^{n}$ and
$\partial E := \overline{E} \setminus E$ denote its frontier. Then
$$ \mathrm{alg.dim}\, (\partial E) < \mathrm{alg.dim}\, (E). $$
\end{proposition}

It is not difficult to strengthen the former proposition as
follows.

\begin{corollary}\label{part1}
Every $\mathcal{L}$-definable set is a finite disjoint union of
locally closed sets.
\end{corollary}

Quantifier elimination due to Pas~\cite[Theorem~4.1]{Pa1} (and
more precisely, elimination of $K$-quantifiers) enables
translation of the language of Denef--Pas on $K$ into a language
described above equipped with the topological system wherein
$\tau_{n}$ is the $K$-topology, $n \in \mathbb{N}$. Indeed, we
must augment the language of valued rings by adding extra relation
symbols for the inverse images under the valuation and angular
component map of relations on the value group and residue field.
More precisely, we must add the names of sets of the form
$$ \{ a \in K^{n}: (v(a_{1},\ldots,v(a_{n}))) \in P \} $$
and
$$ \{ a \in K^{n}: (\overline{ac}\, a_{1},\ldots,\overline{ac}\,
   a_{n}) \in Q \}, $$
where $P$ and $Q$ are definable subsets (in the auxiliary sorts of
the Denef--Pas language) of $\Gamma^{n}$ and $\Bbbk^{n}$,
respectively.

\vspace{1ex}

Summing up, the foregoing results apply in the case of Henselian
non-trivially valued fields with the three-sorted language of
Denef--Pas.

\section{Proof of the main theorem}

Consider an $\mathcal{L}^{P}$-definable function $f: A \to
\mathbb{P}^{1}(K)$ and its graph
$$ E := \{ (x,f(x)): x \in A \} \subset K^{n} \times \mathbb{P}^{1}(K). $$
We shall proceed with induction with respect to the dimension
$$ d = \dim A = \dim \, E $$
of the source and graph of $f$. By Corollary~\ref{part1}, we can
assume that the graph $E$ is a locally closed subset of $K^{n}
\times \mathbb{P}^{1}(K)$ of dimension $d$ and that the conclusion
of the theorem holds for functions with source and graph of
dimension $< d$.

\vspace{1ex}

Let $F$ be the closure of $E$ in $K^{n} \times \mathbb{P}^{1}(K)$
and $\partial E := F \setminus E$ be the frontier of $E$. Since
$E$ is locally closed, the frontier $\partial E$ is a closed
subset of $K^{n} \times \mathbb{P}^{1}(K)$ as well. Let
$$ \pi: K^{n} \times \mathbb{P}^{1}(K) \longrightarrow K^{n} $$
be the canonical projection. Then, by virtue of the closedness
theorem (\cite[Theorem~3.1]{Now}), the images $\pi(F)$ and
$\pi(\partial E)$ are closed subsets of $K^{n}$. Further,
$$ \dim \, F = \dim \, \pi(F) = d $$
and
$$ \dim \, \pi(\partial E) \leq \dim \, \partial E < d; $$
the last inequality holds by Proposition~\ref{front-eq}. Putting
$$ B := \pi(F) \setminus \pi(\partial E) \subset \pi(E) = A, $$
we thus get
$$ \dim \, B = d \ \ \text{and} \ \ \dim \, (A \setminus B) < d.
$$
Clearly, the set
$$ E_{0} := E \cap (B \times \mathbb{P}^{1}(K)) = F \cap (B \times
   \mathbb{P}^{1}(K)) $$
is a closed subset of $B \times \mathbb{P}^{1}(K)$ and is the
graph of the restriction
$$ f_{0}: B \longrightarrow \mathbb{P}^{1}(K) $$
of $f$ to $B$. Again, it follows immediately from the closedness
theorem that the restriction
$$ \pi_{0} : E_{0} \longrightarrow B $$
of the projection $\pi$ to $E_{0}$ is a definably closed map.
Therefore $f_{0}$ is a continuous function. But, by the induction
hypothesis, the restriction of $f$ to $A \setminus B$ satisfies
the conclusion of the theorem, whence so does the function $f$.
This completes the proof of Theorem~\ref{piece}.
\hspace*{\fill}$\Box$

\vspace{2ex}

Finally, let us mention that every $\mathcal{L}^{P}$-definable
continuous function $f: A \to K$ on a closed bounded subset $A$ of
$K^{n}$ is H\"{o}lder continuous, as proven in our recent
paper~\cite{Now1}.

\vspace{2ex}

\begin{small}
Institute of Mathematics

Faculty of Mathematics and Computer Science

Jagiellonian University


ul.~Profesora \L{}ojasiewicza 6

30-348 Krak\'{o}w, Poland

{\em e-mail address: nowak@im.uj.edu.pl}
\end{small}

\end{document}